\documentclass[reqno,11pt]{amsart}

\usepackage{amssymb}
\newtheorem{Proposition}{Proposition}[section]
\newtheorem{Definition}[Proposition]{Definition}
\newtheorem{Lemma}[Proposition]{Lemma}
\newtheorem{Theorem}[Proposition]{Theorem}

\DeclareMathOperator{\Val}{Val}
\DeclareMathOperator{\vol}{vol}
\DeclareMathOperator{\sgn}{sgn}
\DeclareMathOperator{\Alt}{Alt}
\DeclareMathOperator{\Gr}{Gr}

\DeclareMathOperator{\kl}{Kl}
\DeclareMathOperator{\nc}{nc}

\title{A Hadwiger-type theorem for the special unitary group}
\author{Andreas Bernig}
 
\email{andreas.bernig@unifr.ch}
\address{D\'epartement de Math\'ematiques, Chemin du Mus\'ee 23, 1700 Fribourg, Switzerland}

\begin{document}

\begin{abstract}
The dimension of the space of $SU(n)$ and translation invariant continuous valuations on $\mathbb{C}^n, n \geq 2$ is computed. For even $n$, this dimension equals $(n^2+3n+10)/2$; for odd $n$ it equals $(n^2+3n+6)/2$. An explicit geometric basis of this space is constructed. The kinematic formulas for $SU(n)$ are obtained as corollaries. 
\end{abstract}

\thanks{
\\ Supported
  by the Schweizerischer Nationalfonds grants SNF PP002-114715/1 and 200020-113199.}
\maketitle 
%-------------------------------------------------------------------------
\section{Introduction and statement of theorems}

Let $V$ be a finite-dimensional vector space and denote by $\mathcal{K}(V)$ the set of compact convex subsets of $V$. A valuation on $V$ is a map $\mu: \mathcal{K}(V) \to \mathbb{C}$ which is finitely additive in the following sense: 
\begin{displaymath}
 \mu(K \cup L)+\mu(K \cap L)=\mu(K)+\mu(L)
\end{displaymath}
whenever $K, L, K \cup L \in \mathcal{K}(V)$. 

A valuation $\mu$ is continuous if it is continuous with respect to the Hausdorff topology on $\mathcal{K}(V)$. It is translation invariant if $\mu(v+K)=\mu(K)$ for all $v \in V$. 

Now suppose $V$ is a Euclidean vector space. Then $\mu$ is called motion invariant if $\mu(\bar g K)=\mu(K)$ for all Euclidean motions $\bar g$.  
 
Hadwiger's famous characterization theorem states that the space $\Val^{SO(n)}$ of motion invariant continuous valuations is of dimension $n+1$, where $n=\dim V$. The only natural choice (up to scale) of a basis of $\Val^{SO(n)}$ consists of the intrinsic volumes $\mu_0,\ldots,\mu_n$. From Hadwiger's theorem, the array of kinematic formulas, mean projection formulas, additive kinematic formulas and many other results can be obtained in an elegant and simple way. 

It is natural to weaken the hypotheses of motion invariance. The theory of translation invariant continuous valuations is a very rich one, see for instance \cite{mcmu77}, \cite{ale00}, \cite{ale01},  \cite{kl00}, \cite{befu06}, \cite{ale07}. The space of translation invariant continuous valuations is an infinite-dimensional Fr\'echet space. 

In view of potential applications to integral geometry, it is natural to consider translation invariant continuous valuations which are invariant under some compact subgroup $G$ of $O(n)$. Let $\Val^G$ denote the corresponding space. Alesker has shown that $\Val^G$ is finite-dimensional if and only if $G$ acts transitively on the unit sphere \cite{ale00}, \cite{ale06}. This is a strong condition on $G$; in fact Borel \cite{bor49}, \cite{bor50} and Montgomery-Samelson \cite{mosa43} gave a complete classification of connected compact Lie groups acting transitively and effectively on the sphere (compare also \cite{bes87}, 7.13). There are six infinite series 
\begin{equation} \label{eq_series}
 SO(n), U(n), SU(n), Sp(n), Sp(n) \cdot U(1), Sp(n) \cdot Sp(1) 
\end{equation}
and three exceptional groups 
\begin{equation} \label{eq_exceptions}
 G_2, Spin(7),Spin(9). 
\end{equation}

The computation of the dimension of $\Val^G$ and the determination of the kinematic formulas for a group $G$ from this list is a very important question in modern integral geometry. Let us mention what is known in this context. 

Hadwiger's theorem solves the case $G=SO(n)$. The kinematic formulas for this group (which were first proved by Blaschke, Chern and Santal\'o) are corollaries to this theorem. 

The next interesting case is $G=U(n)$ acting on $\mathbb{C}^n$. For $n \leq 3$, Park \cite{pa02} computed the dimension of $\Val^{U(n)}$ and the kinematic formulas. Alesker \cite{ale01} gave the result in the general case: 
\begin{displaymath}
\dim \Val^{U(n)}=\binom{n+2}{2}.
\end{displaymath}
Two different natural and geometric bases of $\Val^{U(n)}$ were constructed in \cite{ale03b}. 

In contrast to the case of the orthogonal group, the computation of the kinematic formulas requires a lot more work and has been completed only very recently \cite{fu06}, \cite{befu08}. 

Alesker \cite{ale04} showed that $\dim \Val^{SU(2)}=10$. The kinematic formulas for $SU(2)$ have been obtained in \cite{be06}. 

For the groups $Sp(n), Sp(n) \cdot U(1), Sp(n) \cdot Sp(1) $ as well as for the exceptional groups $G_2, Spin(7),Spin(9)$, a Hadwiger-type theorem is still unknown. Some partial results have been obtained by Alesker: several $Sp(n)\cdot Sp(1)$-invariant valuations on $\mathbb{H}^n$ were constructed in \cite{ale05}; a new  $Spin(9)$-invariant valuation on $\mathbb{R}^{16}$ is given in \cite{ale07b}. 
 
The present paper is devoted to the study of the integral geometry of $SU(n)$ for all $n \geq 2$. We compute the dimension of the space $\Val^{SU(n)}$ and derive all kinematic formulas.

Before stating our main results, we need some notation. 

Recall that a translation invariant valuation $\mu$ on a vector space $V$ is said to be of degree $k$ if $\mu(tK)=t^k \mu(K)$ for all $t \geq 0$ and all compact convex sets $K$. Each continuous translation invariant valuation on $V$ can be uniquely decomposed as a sum of valuations of degrees $0,1,\ldots,\dim V$ \cite{mcmu77}. We let $\Val_k(V)$ denote the space of continuous translation invariant valuations of degree $k$.  

\begin{Definition} \label{def_weight}
A valuation $\mu \in \Val^{SU(n)}$ has {\bf weight} $l$ if 
\begin{displaymath}
\mu(gK)=\det(g)^l \mu(K)
\end{displaymath}
for all $g \in U(n)$. We write $\Val^{SU(n),l}$ for the space of valuations of weight $l$. Clearly $\Val^{SU(n),0}=\Val^{U(n)}$. 
\end{Definition}

Taking $g=-1$, we see that a valuation $\mu$ of weight $l$ is even (i.e. $\mu(-K)=\mu(K)$ for all $K$) if $nl$ is even and odd (i.e. $\mu(-K)=-\mu(K)$ for all $K$) if $nl$ is odd. 

\begin{Proposition} \label{prop_splitting}
\begin{enumerate}
 \item The space $\Val^{SU(n)}$ admits a splitting
\begin{displaymath}
\Val^{SU(n)}=\bigoplus_{l=-2}^2 \Val^{SU(n),l}.
\end{displaymath}
In particular, the weight of a valuation can only be $0, \pm 1, \pm 2$. 
\item The weight of an Alesker product is the sum of the weights of the factors. 
\item If $\Phi:\Val^{SU(n)} \to \Val^{SU(n)}$ is any linear operator commuting with the action of $U(n)$, then $\Phi$ preserves the weight.
\item There are natural isomorphisms, given by complex conjugation, 
\begin{equation} \label{eq_weight_isom}
\Val^{SU(n),1} \cong \Val^{SU(n),-1}, \quad \Val^{SU(n),2} \cong \Val^{SU(n),-2}.
\end{equation}
\end{enumerate}
\end{Proposition}

Let $W \subset \mathbb{C}^n$ be a (real) subspace of dimension $n$. By work of Tasaki, one can associate $m=\lfloor n/2\rfloor$ K\"ahler angles $0 \leq \theta_1 \leq \ldots \leq \theta_m \leq \frac{\pi}{2}$ to $W$. These numbers describe the orbits of the $U(n)$-action on $\Gr_n(\mathbb{C}^n)$. We refer to Section \ref{sec_theta} and \cite{tas00} for their definition. 

We choose an orthonormal basis $w_1,w_2,\ldots,w_n$ of $W$. If the restriction of the symplectic form of $\mathbb{C}^n$ to $W$ is non-degenerated (this is the case if and only if $n$ is even and all K\"ahler angles are strictly less than $\frac{\pi}{2}$), we want $w_1 \wedge \ldots \wedge w_n$ to induce the same orientation as the symplectic form. 
\begin{Definition}
The {\bf $\Theta$-invariant} of $W$ is the number 
\begin{displaymath}
 \Theta(W):=\det(w_1,\ldots,w_n)
\end{displaymath}
which is a complex number if the restriction of the symplectic form to $W$ is non-degenerated and which is an element of $\mathbb{C}/\{\pm 1\}$ otherwise. 
\end{Definition}

\begin{Proposition} \label{prop_orbitspace}
\begin{enumerate}
\item Let $W_1,W_2$ be two $n$-dimensional subspaces of $\mathbb{C}^n$. There exists $g \in SU(n)$ with $gW_1=W_2$ if and only if $W_1$ and $W_2$ have the same K\"ahler angles and the same $\Theta$-invariant. 
\item If $W \in \Gr_n(\mathbb{C}^n)$ has K\"ahler angles $\theta_1,\ldots,\theta_m$ and $\Theta:=\Theta(W)$ then 
\begin{equation} \label{eq_norm_Theta}
 |\Theta|=\prod_{j=1}^m \sin(\theta_j).
\end{equation}
Conversely, given $0 \leq \theta_1 \leq \ldots \leq \theta_m \leq \pi/2$ and $\Theta$ satisfying \eqref{eq_norm_Theta}, there exists $W \in \Gr_n(\mathbb{C}^n)$ with $\theta_j(W)=\theta_j, \Theta(W)=\Theta$.  
\item For all $W \in \Gr_n(\mathbb{C}^n)$
\begin{equation} \label{eq_Theta_perp}
\Theta(W^\perp)=\Theta(W).
\end{equation} 
\item If $k \neq n$, then the orbits of $U(n)$ and $SU(n)$ on $\Gr_k(\mathbb{C}^n)$ are the same. 
\end{enumerate}
\end{Proposition}

An even valuation $\mu \in \Val_k(\mathbb{C}^n)$ is completely determined by its restriction to $k$-dimensional subspaces. More precisely, if $W \in \Gr_k(\mathbb{C}^n)$, then $\mu|_W$ is a multiple of the $k$-dimensional volume on $W$. The proportionality factor is denoted by $\kl_\mu(W)$ and $\kl_\mu \in C(\Gr_k(\mathbb{C}^n))$ is called {\bf Klain function} of $\mu$. The resulting map $\kl:\Val_k^+(\mathbb{C}^n) \to C(\Gr_k(\mathbb{C}^n)), \mu \mapsto \kl_\mu$ is injective by a result of Klain \cite{kl00}. 

In general, it is rather difficult to write down a valuation with given Klain function. However, a {\bf constant coefficient valuation} (see \cite{befu08} and Subsection \ref{subsection_construction_phi2} for the definition) can be easily recovered from its Klain function. If $\mu$ is such a valuation and $P \subset \mathbb{C}^n$ a convex polytope, then 
\begin{displaymath}
 \mu(P)=\sum_{F, \dim F=k} \gamma(F) \vol(F) \kl_\mu(W_F),
\end{displaymath}
where $F$ runs over all $k$-dimensional faces of $P$; $W_F \in \Gr_k(\mathbb{C}^n)$ is the linear space parallel to $F$ and $\gamma(F)$ is the normalized volume of the exterior angle at $F$. 

Our main theorem is the following structure theorem for $\Val^{SU(n)}$. 

\begin{Theorem} \label{thm_hadwiger}
There exists a unique constant coefficient valuation $\phi_2 \in \Val_n^{SU(n)}$ with Klain function  
\begin{equation} \label{eq_klain_phi2}
\kl_{\phi_2} =\Theta^2.
\end{equation}
$\phi_2$ is even, of degree $n$ and spans $\Val^{SU(n),2}$. 
If $n=2m$ is even, there exists a unique constant coefficient valuation $\phi_1 \in \Val_n^{SU(n)}$ with Klain function  
\begin{equation} \label{eq_klain_phi1}
\kl_{\phi_1}(W) =\Theta(W)  \prod_{j=1}^m \cos(\theta_j(W)), \quad W \in \Gr_n(\mathbb{C}^n).
\end{equation}
$\phi_1$ is even, of degree $n$ and spans $\Val^{SU(n),1}$. If $n$ is odd, then $\Val^{SU(n),1}=0$.  
In particular, there are no odd invariant valuations and 
\begin{align*}
 \Val_k^{SU(n)} & =\Val_k^{U(n)} \quad \text{ if } k \neq n;\\
\dim \Val_n^{SU(n)} & = \dim \Val_n^{U(n)}+4 \quad \text{ if } n \equiv 0 \mod 2;\\
\dim \Val_n^{SU(n)} & = \dim \Val_n^{U(n)}+2 \quad \text{ if } n \equiv 1 \mod 2.
\end{align*}
\end{Theorem}

It is not clear a priori why $SU(n)$-invariant valuations of degree $k \neq n$ are $U(n)$-invariant. Also the fact that such valuations are even is not trivial if $n$ is odd. We do not know if there is a geometric proof of these facts. Note, however, that the second statement implies the first one: since the $U(n)$-orbit and the $SU(n)$-orbit on $\Gr_k(\mathbb{C}^n)$ agree for $k \neq n$ by Proposition \ref{prop_orbitspace}, Klain's injectivity theorem implies that {\it even} $SU(n)$-invariant valuations of degree $k \neq n$ are $U(n)$-invariant. 

Our approach is based on the fact that $SU(n)$-invariant valuations can be identified with a quotient of the space of $SU(n)$-invariant differential forms on the unit sphere bundle of $\mathbb{C}^n$. This is a consequence of Alesker's irreducibility theorem \cite{ale01} and the kernel theorem of \cite{bebr07}. Here it is important that $SU(n)$ acts transitively on the unit sphere, compare \cite{ale06} and \cite{fu04} for more information.     

Our main application is to the integral geometry of the group $SU(n)$ acting on $\mathbb{C}^n$. Recall that for any group $G$ from the lists \eqref{eq_series} and \eqref{eq_exceptions}, one may define an injection $k_G:\Val^G \to \Val^G \otimes \Val^G$ by 
\begin{displaymath}
k_G(\mu)(K,L)=\int_{\bar G} \mu(K \cap \bar g L)d \bar g, \quad K,L \in \mathcal{K}(V). 
\end{displaymath}
Then $k_G$ is a a cocommutative, coassociative coproduct. 

A detailed study of $k_{U(n)}$ (generalizing results of Park \cite{pa02}, Tasaki \cite{tas00}, Alesker \cite{ale03b} and Fu \cite{fu06}) is contained in \cite{befu08}.   
In the next theorem, we identify $\Val^{U(n)}$ with the corresponding subspace in $\Val^{SU(n)}$. We let $\omega_n$ be the volume of the $n$-dimensional unit ball. 

\begin{Theorem} \label{thm_kinematic}
The principal kinematic formulas for $U(n)$ and $SU(n)$ are related by 
\begin{displaymath}
k_{SU(n)}(\chi)  =k_{U(n)}(\chi)+\frac{1}{2^n}(\phi_1 \otimes \bar \phi_1+\bar \phi_1 \otimes \phi_1)+\frac{\omega_n^2}{(n+2) 2^{2n-1} \omega_{2n}} (\phi_2 \otimes \bar \phi_2+\bar \phi_2 \otimes \phi_2) 
\end{displaymath}
if $n$ is even and by 
\begin{displaymath}
k_{SU(n)}(\chi)  =k_{U(n)}(\chi)-\frac{\omega_n^2}{(n+2) 2^{2n-1} \omega_{2n}}  (\phi_2 \otimes \bar \phi_2+\bar \phi_2 \otimes \phi_2) 
\end{displaymath}
if $n$ is odd. If $\mu$ is a $U(n)$-invariant valuation of degree $k>0$, then $k_{SU(n)}(\mu)=k_{U(n)}(\mu)$. If $n$ is even, then 
\begin{displaymath}
k_{SU(n)}(\phi_1) = \phi_1 \otimes \vol+\vol \otimes \phi_1. 
\end{displaymath}
For all $n \geq 2$
\begin{displaymath}
k_{SU(n)}(\phi_2) = \phi_2 \otimes \vol+\vol \otimes \phi_2.
\end{displaymath}
\end{Theorem}

The interested reader may rewrite this theorem in terms of the algebra structure of $\Val^{SU(n)}$. By \cite{fu06} and Theorem \ref{thm_hadwiger}, this graded algebra is generated by a unitarily invariant valuation $t$ of degree $1$, a unitarily invariant valuation $s$ of degree $2$, the two $SU(n)$-invariant valuations $\phi_2,\bar \phi_2$ of degree $n$ and, if $n$ is even, the two $SU(n)$-invariant valuations $\phi_1,\bar \phi_1$ of degree $n$. The relations between these elements are computed in Section \ref{sec_kinematic}.

\subsection*{Plan of the paper}
In Section \ref{sec_theta} we introduce the $\Theta$-invariant and prove Proposition \ref{prop_orbitspace}. In Section \ref{sec_forms} we construct a generating set of the algebra of $SU(n)$-invariant differential forms on $S\mathbb{C}^n=\mathbb{C}^n \times S^{2n-1}$. We also establish some relations between such forms. The following Section \ref{sec_weight} contains the proof of Proposition \ref{prop_splitting}. The main part of the paper is the proof of Theorem \ref{thm_hadwiger} in Section \ref{sec_hadwiger}. In the final Section \ref{sec_kinematic}, we recall the definition of Alesker's product of valuations and its relation with kinematic formulas, we prove Theorem \ref{thm_kinematic} and we derive an additive kinematic formula. 

\subsection*{Acknowledgements} I wish to thank Semyon Alesker and Joseph Fu for very useful remarks on a first version of this manuscript. Proposition \ref{prop_splitting} (c) and some other improvements were suggested by the anonymous referee.

%----------------------------------------------------------------

\section{The orbit space of $\Gr_k(\mathbb{C}^n)$ under the $SU(n)$-action}
\label{sec_theta}

By work of Tasaki \cite{tas00} each $U(n)$-orbit in $\Gr_n(\mathbb{C}^n)$ is described by $m:=\lfloor n/2\rfloor$ K\"ahler angles $0 \leq \theta_1 \leq \ldots \leq \theta_m \leq \pi$. 

Given $W \in \Gr_n(\mathbb{C}^n)$, we consider the composition $J_W:= \pi_W\circ J|_W:W \to W$ of the orthogonal projection $\pi_W$ with the multiplication $J$ by $\sqrt{-1}$. Since $J_W$ is
skew-symmetric and $\|J_W\| \leq 1$, each eigen-value of  $J_W$ is purely
imaginary and has absolute value at most $1$. The multiple K\"ahler angle of  $W$ is
the $m$-tuple $(\theta_1,\ldots,\theta_m)$ with
\begin{displaymath}
0\le\theta_1\le \theta_2\le \dots\le \theta_m\le \frac \pi 2
\end{displaymath}
such that 
$\left\{\pm \cos(\theta_1) i ,\ldots,\pm \cos(\theta_m) i\right\}$
(and $0$ if $n$ is odd) are the eigenvalues of $\pi_W\circ J|W$.

Two $n$-dimensional spaces belong to the same $U(n)$-orbit if and only if their multiple K\"ahler angles agree \cite{tas00}. Fixing a hermitian basis $e_1,\ldots,e_n$ of $\mathbb{C}^n$, a typical space in such an orbit is given by 
\begin{equation} \label{eq_un_orbit}
W:=\bigoplus_{j=1}^m \left[\mathbb{R}e_{2j-1} \oplus \mathbb{R}(\cos(\theta_j) ie_{2j-1}+\sin(\theta_j)e_{2j}) \right]\underbrace{\left(\bigoplus \mathbb{R}e_n\right)}_{\text{only if } n \equiv 1 \mod 2}.
\end{equation}

Now we turn our attention to the action of $SU(n)$. 

Let $W \in \Gr_n(\mathbb{C}^n)$ have K\"ahler angles $\theta_1,\ldots,\theta_m$. If $\theta_m < \pi/2$ and $n$ is even, then the restriction of the K\"ahler form is non-degenerated and its $m$-th power defines an orientation of $W$. In all other cases, we fix an arbitrary orientation of $W$. 

Fix a positively oriented orthogonal basis (w.r.t. to the Euclidean structure on $\mathbb{C}^n$) $w_1,\ldots,w_n$ of $W$ and define 
\begin{displaymath}
 \Theta(W):=\det(w_1,\ldots,w_n) \in \mathbb{C}.
\end{displaymath}

If $\theta_m < \pi/2$ and $n$ is even, $\Theta(W) \in \mathbb{C}$ is independent of the choice of the positively oriented orthonormal basis. In the other cases, $\Theta(W) \in \mathbb{C}/\{\pm 1\}$ is independent of the choice of the orientation and of the choice of the orthonormal basis. 

We call $\Theta(W)$ the {\bf $\Theta$-invariant} of $W$. 

Remark: As was pointed out to us by J. Fu, the restriction of $\Theta^2$ to the Lagrangian Grassmannian has values in $S^1$ and is a primitive of the Maslov $1$-form (compare \cite{mcduffsal95}, page 53). 

\proof[Proof of Proposition \ref{prop_orbitspace}]
Fix a positively oriented orthonormal basis $w_1,\ldots,w_n$ of $W \in \Gr_n(\mathbb{C}^n)$. If $g \in U(n)$, then $gw_1,\ldots,gw_n$ is a positively oriented basis of $gW$. Hence 
\begin{equation} \label{eq_Theta}
\Theta(gW)=\det(g) \Theta(W).
\end{equation} 
In particular, $\Theta$ is $SU(n)$-invariant. 

Conversely, suppose that $W_1$ and $W_2$ in $\Gr_n(\mathbb{C}^n)$ have the same K\"ahler angles and that $\Theta(W_1)=\Theta(W_2)$. 

If $\Theta(W_1)=0$, then the complex subspace $W_1^\mathbb{C}$ generated by $W_1$ is of positive codimension. Let $g \in U(n)$ be such that $gW_1=W_2$. Let $g_0$ be the element which acts by multiplication by $\det(g)^{-1}$ on a one-dimensional complex subspace in the complement of $W_1^\mathbb{C}$ and identically elsewhere. Then $g_0$ fixes $W_1$. Hence $gg_0W_1=W_2$ and $gg_0 \in SU(n)$. 

Let us now assume that $\Theta(W_1) \neq 0$. If $n$ is even and $\theta_m < \pi/2$, then \eqref{eq_Theta} shows that $g \in SU(n)$. 

If $n$ is odd or if $\theta_m =\pi/2$, we can only deduce that $\det(g) \in \{\pm 1\}$. In both cases, there exists $w_1 \in W$ with $Jw_1 \perp W$. We complete $w_1$ to an orthonormal basis $\{w_1,\ldots,w_n\}$ of $W$. Then $w_2,\ldots,w_n$ belong to the orthogonal complement of the complex subspace $\mathbb{C}w_1$. Let $g_0 \in U(n)$ be the element which acts by multiplication by $\det(g)$ on $\mathbb{C}w_1$ and by the identity on $(\mathbb{C}w_1)^\perp$. Then $g_0$ fixes $W_1$. It follows that $gg_0W_1=W_2$ and $gg_0 \in SU(n)$. This proves the first part. 

For the second assertion, we consider the space $W$ defined in \eqref{eq_un_orbit}. Clearly 
\begin{displaymath} 
\Theta(W)=\prod_{j=1}^m \sin(\theta_j).
\end{displaymath}
If $W_1$ is any other space with these multiple K\"ahler angles, there exists $g \in U(n)$ with $gW=W_1$. Then 
\begin{displaymath} 
|\Theta(W_1)|=|\det(g)\Theta(W)|=|\Theta(W)|=\prod_{j=1}^m \sin(\theta_j). 
\end{displaymath}

Conversely, given any $\Theta$ with $|\Theta|=\prod_{j=1}^m \sin(\theta_j)$, we find $\xi \in S^1$ with $\xi \cdot  \prod_{j=1}^m \sin(\theta_j)=\Theta$. Taking $g \in U(n)$ of determinant $\xi$ and setting $W_1:=gW$, we obtain that $\Theta(W_1)=\Theta$. 

It suffices to show the third assertion for one space in each $U(n)$-orbit. As before, we take the space $W$ of \eqref{eq_un_orbit}. The $\Theta$-invariant of $W$ equals $\prod_{j=1}^m \sin(\theta_j)$. The $\Theta$-invariant of its orthogonal complement 
\begin{displaymath}
W^\perp=\bigoplus_{j=1}^m \left[\mathbb{R}(-\sin(\theta_j) i e_{2j-1}+\cos(\theta_j)e_{2j}) \oplus \mathbb{R} i e_{2j}\right] \underbrace{\left(\bigoplus \mathbb{R} ie_n\right)}_{\text{only if } n \equiv 1 \mod 2};
\end{displaymath}
is again $\prod_{j=1}^m \sin(\theta_j)$, which proves the third assertion. 

If $W \in \Gr_k(\mathbb{C}^n)$ for $k<n$, then the complex space generated by $W$ has positive codimension. It follows that for any $u \in S^1$ there is an element $g \in U(n)$ which stabilizes $W$ and which has determinant $u$. Therefore the $SU(n)$-orbit and the $U(n)$-orbit on $\Gr_k(\mathbb{C}^n)$ are the same. By passing to the orthogonal complements, the same is true for $k>n$. 

\endproof
%----------------------------------------------------------------

\section{Invariant differential forms}
\label{sec_forms}

Let $S\mathbb{C}^n=\mathbb{C}^n \times S^{2n-1}$ be the sphere bundle of $\mathbb{C}^n$. We consider the diagonal action of $U(n)$ on $S\mathbb{C}^n$. Park \cite{pa02} described the algebra of $U(n)$-invariant forms on $S\mathbb{C}^n$ as follows. 

Let us introduce complex coordinates $z_j=x_j+iy_j, j=1,\ldots,n$ on $\mathbb{C}^n$ and use induced coordinates $(z_j,\zeta_j=\xi_j+i \eta_j)$ on $\mathbb{C}^n \times \mathbb{C}^n$. We consider $S\mathbb{C}^n$ as a subset of $\mathbb{C}^n \times \mathbb{C}^n$ and define real-valued differential forms by  
\begin{align*}
\alpha & = \frac12 \sum_{j=1}^n (\zeta_j d\bar z_j+\bar \zeta_j dz_j) \\
\beta & = \frac{i}{2} \sum_{j=1}^n (\zeta_j d\bar z_j-\bar \zeta_j dz_j)\\
\gamma & = \frac{i}{2} \sum_{j=1}^n (\zeta_j d\bar \zeta_j-\bar \zeta_j d \zeta_j) \\
\theta_0 & = \frac{i}{2} \sum_{j=1}^n d\zeta_j \wedge d\bar\zeta_j \\
\theta_s- i \theta_1 & = \sum_{j=1}^n dz_j \wedge d\bar \zeta_j-\beta \wedge \gamma+i \alpha \wedge \gamma \\
\theta_2 & := \frac{i}{2} \sum_{j=1}^n dz_j \wedge d\bar z_j-\alpha \wedge \beta.
\end{align*}

For all computations not involving exterior differential, it is convenient to use the fact that $SU(n)$ acts transitively on the unit sphere and that an invariant form is thus determined by its value at the special point $(0,e_1) \in S\mathbb{C}^n$. At this point, the above forms are given by 
\begin{align*}
\alpha & = dx_1\\
\beta & = dy_1 \\
\gamma & = d\eta_1\\
\theta_0 & = \sum_{j=2}^n d\xi_j \wedge d\eta_j\\
\theta_s & =\sum_{j=2}^n dx_j \wedge d\xi_j+dy_j \wedge d\eta_j\\
\theta_1 & =\sum_{j=2}^n dx_j \wedge d\eta_j-dy_j \wedge d\xi_j\\
\theta_2 & = \sum_{j=2}^n dx_j \wedge dy_j.
\end{align*}

\begin{Proposition} \cite{pa02}\\ \label{prop_park_classification}
The algebra of $U(n)$-invariant forms on $S\mathbb{C}^n$ is generated by $\alpha,\beta,\gamma,\theta_0, \theta_1,\theta_2$ and $\theta_s$. They satisfy
\begin{align*}
d\alpha & =-\beta \wedge \gamma-\theta_s\\
d\beta & = \alpha \wedge \gamma+\theta_1\\
d\gamma & = 2\theta_0\\
d\theta_0 & = 0\\
d\theta_1 & =2\alpha \wedge \theta_0+ \gamma \wedge \theta_s\\
d\theta_2 & = \alpha \wedge \theta_1+\beta \wedge \theta_s\\
d\theta_s & = 2 \beta \wedge \theta_0-\gamma \wedge \theta_1. 
\end{align*}
\end{Proposition}

Let us introduce $n-1$-forms $\chi_k, k=0,\ldots,n-1$ by setting 
\begin{displaymath}
\chi_k=\frac{1}{k!(n-k-1)!} \sum_{\pi \in \mathcal{S}_n} \sgn(\pi) \zeta_{\pi(1)} dz_{\pi(2)} \wedge \ldots \wedge dz_{\pi(k+1)} \wedge d\zeta_{\pi(k+2)} \wedge \ldots \wedge d\zeta_{\pi(n)}.
\end{displaymath}
Here $\mathcal{S}_n$ is the group of permutations of $\{1,\ldots,n\}$ and $\sgn(\pi)$ denotes the sign of a permutation $\pi$. For convenience, we set $\chi_k:=0$ if $k \not\in \{0,\ldots,n-1\}$. 

\begin{Proposition} \label{prop_invariant_forms}
The algebra of (complex-valued) $SU(n)$-invariant forms on $S\mathbb{C}^n$ is generated by $\alpha,\beta,\gamma,\theta_0, \theta_1,\theta_2, \theta_s$, $\chi_k$ and $\bar \chi_k$ ($k=0,\ldots,n-1$). Moreover,
\begin{equation} \label{eq_differential_chi}
 d\chi_k = (n-k) ((\alpha+i\beta) \wedge \chi_{k-1}+i \gamma \wedge \chi_k). 
\end{equation}
\end{Proposition}

For the proof, we need the following lemma. 
\begin{Lemma}
Consider $V:=\mathbb{C}^n \oplus \mathbb{C}^n$ with the diagonal action of $SU(n)$. Let $e_1^0,\ldots,e_n^0,e_1^1,\ldots,e_n^1$ be the standard basis of $V$, and denote the dual basis by $e_1^{0*},\ldots,e_n^{0*},e_1^{1*},\ldots,e_n^{1*} \in V^*$. Then $\Lambda^*(V^*)^{SU(n)}$ is generated, as an algebra, by the following elements: 
\begin{align*}
\Theta_0 & := \frac{i}{2} \sum_{j=1}^n e^{1*}_j \wedge \bar e^{1*}_j\\
\Theta_2 & :=\frac{i}{2} \sum_{j=1}^n e^{0*}_j \wedge \bar e^{0*}_j\\
\Theta_s & :=\frac12 \sum_{j=1}^n (e^{0*}_j \wedge \bar e^{1*}_j+\bar e^{0*}_j \wedge e^{1*}_j)\\
\Theta_1 & := \frac{i}{2} \sum_{j=1}^n (e^{0*}_j \wedge \bar e^{1*}_j-\bar e^{0*}_j \wedge e^{1*}_j)\\
\Xi_k & := \frac{1}{k!(n-k)!} \sum_{\pi \in \mathcal{S}_n} \sgn(\pi) e^{0*}_{\pi(1)} \wedge \ldots e^{0*}_{\pi(k)} \wedge e^{1*}_{\pi(k+1)} \wedge \ldots \wedge e^{1*}_{\pi(n)}, \quad k=0,\ldots,n\\
\bar \Xi_k & := \frac{1}{k!(n-k)!} \sum_{\pi \in \mathcal{S}_n} \sgn(\pi) \bar e^{0*}_{\pi(1)} \wedge \ldots \bar e^{0*}_{\pi(k)} \wedge \bar e^{1*}_{\pi(k+1)} \wedge \ldots \wedge \bar e^{1*}_{\pi(n)}, \quad k=0,\ldots,n.
\end{align*}
\end{Lemma}

\proof
We adapt the proof of \cite{pa02}, Thm. 2.1.2 to our situation. Let $\phi \in \Lambda^k (V^*)^{SU(n)}$. Write $v=(v^0,v^1)$ for an element $v \in V$. By the first fundamental theorem for $SU(n)$ (\cite{spiv5}, Thm. 64), $\phi(v_1,\ldots,v_k)$ may be written as a polynomial in the real and imaginary parts of the hermitian scalar products 
\begin{displaymath}
(v_i^a,v_j^b), \quad i,j=1,\ldots,k; a,b=0,1;
\end{displaymath}
the determinant functions
\begin{displaymath}
\det(v_{i_1}^{a_1},\ldots,v_{i_n}^{a_n}), \quad i_1,\ldots,i_n \in \{1,\ldots,k\}; a_1,\ldots,a_n \in \{0,1\}
\end{displaymath} 
and their conjugates. Since $\phi$ is multilinear, each lower index $i \in \{1,\ldots,k\}$ appears exactly once in each monomial. 

We apply the alternation operator 
\begin{displaymath}
\Alt: \otimes{}^* V^* \to \Lambda^* V^*
\end{displaymath}
to such an expression. Since $\Alt$ is linear, $\Alt(\psi \otimes \vartheta)=c \Alt(\psi) \wedge \Alt(\vartheta)$ and $\Alt \phi=\phi$, we deduce that $\phi$ is the wedge product of the alternations of the real and imaginary parts of the hermitian scalar products and the determinant functions and their conjugates. These alternations are - up to constants - the elements $\Theta_0,\Theta_1,\Theta_2,\Theta_s,\Xi_k,\bar \Xi_k$. 
\endproof

\proof[Proof of Proposition \ref{prop_invariant_forms}]
First note that $SU(n)$ acts transitively on the unit sphere. If $\omega$ is an invariant form on $S\mathbb{C}^n$, then $\omega$ is determined by its value at the point $(0,e_1)$. The tangent space at this point splits as 
\begin{displaymath}
T_{(0,e_1)}S\mathbb{C}^n= \mathbb{C} e_1 \oplus \bigoplus_{j=2}^n \mathbb{C}e_j \oplus \mathbb{R}ie_1 \oplus \bigoplus_{j=2}^n \mathbb{C}e_j.
\end{displaymath}
The action of the stabilizer at $(0,e_1)$ on the first and third factor is trivial and is the diagonal action on the second and fourth factor. Note that $\alpha + i\beta$ is the projection on the first factor and $\gamma$ is the projection on the third factor. The first statement of the proposition thus follows from the lemma. Equation  \eqref{eq_differential_chi} is easily obtained by comparing the restriction to $T_{(0,e_1)}$ of both sides.   
\endproof

\begin{Proposition} \label{prop_chi_theta}
For $k=1,\ldots,n-1$ 
\begin{align}
 \chi_k \wedge \theta_0 & = -\frac{i}{2} \chi_{k-1} \wedge (\theta_s-i\theta_1)\\
\chi_{k-1} \wedge \theta_2 & = \frac{i}{2} \chi_k \wedge (\theta_s+i\theta_1).
\end{align}
If $k+l \geq n$, then
\begin{equation} \label{eq_relation_chi_theta1}
\chi_k \wedge (\theta_s-i \theta_1)^{l}=0.
\end{equation}
If $l > k$, then
\begin{equation} \label{eq_relation_chi_theta2}
\chi_k \wedge (\theta_s+i \theta_1)^{l}=0.
\end{equation}
\end{Proposition}

\proof
Easy computation.
\endproof

%----------------------------------------------------------------

\section{Weight decomposition}
\label{sec_weight}

\proof[Proof of Proposition \ref{prop_splitting}] 

Recall from Definition \ref{def_weight} that a valuation $\mu \in \Val^{SU(n)}$ has weight $l$ if $\mu(gK)=\det(g)^l \mu(K)$ for all $g \in U(n)$ and all compact convex sets $K$. We shall show that each valuation can be decomposed into valuations of weight $\pm 2,\pm 1$ or $0$. 

Let us say that a differential form $\omega \in \Omega^*(S\mathbb{C}^n)^{SU(n)}$ has weight $l$ if $g^*\omega=(\det g)^l \omega$ for all $g \in U(n)$. Clearly, $\alpha,\beta,\gamma,\theta_0,\theta_1,\theta_2,\theta_s$ are of weight zero (since they are invariant under $U(n)$), while $\chi_k$ is of weight $1$ and $\bar \chi_k$ is of weight $-1$. The weight is additive under wedge products and invariant under differentiation. Each invariant form $\omega$ can be decomposed in a unique way as 
\begin{equation} \label{eq_dec_omega}
\omega=\sum_{l=-\infty}^{\infty} \omega_l,
\end{equation} 
where $\omega_l$ is of weight $l$ and where only finitely many $\omega_l$ are non-zero.  

The only $SU(n)$-invariant valuation in degree $2n$ is the volume, which is of weight $0$. If $\mu$ is an invariant valuation of degree $k<2n$, then we represent $\mu$ by $\omega \in \Omega^{2n-1}(S\mathbb{C}^n)^{\overline{SU(n)}}$, i.e. if $\nc(K)$ is the normal cycle of $K$ (compare \cite{zaeh86} for the normal cycle of compact convex sets), then 
\begin{displaymath}
\mu(K)=\nc(K)(\omega). 
\end{displaymath}

The fact that such an expression exists is a consequence of Alesker's irreducibility theorem \cite{ale01}. Note, however, that $\omega$ is not unique in general. The forms inducing the trivial valuation are characterized in \cite{bebr07}. All we need here is that exact forms and vertical forms (i.e. forms vanishing on the contact distribution) induce the zero valuation; which is a triviality since $\nc(K)$ is a Legendrian cycle. 

In the case $n=2$, any product of $3$ forms of weight $1$ vanishes (since there are only two such forms, $\chi_0$ and $\bar \chi_0$, and they are of degree $1$). In the case $n>2$, we have $3 \deg \chi_k=3(n-1) > 2n-1=\deg \omega$, hence $\omega$ can not contain three or more factors of weight $1$. In both cases it follows that the non-zero terms in the decomposition \eqref{eq_dec_omega} can only appear for $l \in \{0,\pm 1, \pm 2\}$. 

Let $\mu_l$ be the valuation represented by $\omega_l$. Clearly, $\mu_l$ has weight $l$ and $\mu=\sum_{l=-2}^2 \mu_l$. The uniqueness of such a decomposition is clear. This proves (a).  

Statement (b) follows from the fact that a linear automorphism on a vector space $V$ induces an algebra automorphism on $\Val(V)$. 

Statement (c) is trivial: $\mu$ is of weight $l$ if and only if $g \mu=(\det g)^{-l}\mu$ for all $g \in U(n)$. If $\Phi$ commutes with the action of $U(n)$, then 
\begin{displaymath}
 g \Phi(\mu)=\Phi(g \mu)=(\det g)^{-l} \Phi(\mu), \quad g \in U(n);
\end{displaymath}
hence $\Phi(\mu)$ is again of weight $l$. 

In order to prove (d), we note that if $\mu$ is of weight $l$, then $\bar \mu$ is of weight $-l$: 
\begin{displaymath}
\bar \mu(gK)=\overline{\det g}^l \bar \mu(K)=(\det g)^{-l} \bar \mu(K). 
\end{displaymath}
\endproof

%----------------------------------------------------------------

\section{Classification of $SU(n)$-invariant valuations}
\label{sec_hadwiger}

In this section, we prove Theorem \ref{thm_hadwiger}.

\subsection{Construction of $\phi_2$} 
\label{subsection_construction_phi2}

Let $\omega_n$ be the volume of the $n$-dimensional unit ball. We claim that the valuation 
\begin{displaymath}
\phi_2(K):=\frac{(-1)^{n+1}i}{n \omega_n} \nc(K)(\beta \wedge \chi_0 \wedge \chi_{n-1})
\end{displaymath}
satisfies \eqref{eq_klain_phi2}. Note that 

\begin{displaymath}
\beta \wedge \chi_0 \wedge \chi_{n-1} \equiv (-1)^n i dz_1 \wedge \ldots dz_n \wedge \sum_{j=1}^n (-1)^{j+1} \zeta_j d\zeta_1 \wedge \ldots \wedge \hat{d\zeta_j} \wedge \ldots \wedge d\zeta_n \mod \alpha.
\end{displaymath}

It follows that $\phi_2$ is a constant coefficient valuation in the sense of \cite{befu08}. In fact, the right hand side of this relation clearly extends to a form on $\mathbb{C}^n \oplus \mathbb{C}^n$ whose exterior differential has constant coefficients. 

Let $W \in \Gr_n(\mathbb{C}^n)$ and denote the unit ball inside $W$ by $B_W$. 

The part of bidegree $(n,n-1)$ of $\nc(B_W)$ is given by $[[B_W]] \times [[\partial B_W]]$. 

It follows that 
\begin{align*}
\phi_2(B_W) & = \frac{(-1)^{n+1}i}{n\omega_n} \nc(B_W)(\beta \wedge \chi_0 \wedge \chi_{n-1})\\
& =\frac{1}{n\omega_n} [[B_W]](dz_1 \wedge \ldots dz_n) [[\partial B_W]]\left(\sum_{j=1}^n (-1)^{j-1} \zeta_j d\zeta_1 \wedge \ldots \wedge \hat{d\zeta_j} \wedge \ldots \wedge d\zeta_n\right) \\
& = \vol(B_W) \Theta(W) \Theta(W^\perp). 
\end{align*}

Using Proposition \ref{prop_orbitspace} we thus obtain 
\begin{displaymath}
 \kl_{\phi_2}(W)=\Theta(W)^2. 
\end{displaymath}

\subsection{Classification of invariant valuations of weight $2$} A valuation of weight $2$ must be represented by a form $\omega$ with $g^*\omega=(\det g)^2 \omega$ for all $g \in U(n)$. The vector space of such forms is generated by $\alpha \wedge \chi_j \wedge \chi_k$ (which induces the zero valuation), $\beta \wedge \chi_j \wedge \chi_k$ and $\gamma \wedge \chi_j \wedge \chi_k$ with $j,k=0,\ldots,n-1$. 

It follows readily from the definition that $\chi_j \wedge \chi_k=0$ unless $j+k=n-1$ and that the forms $\chi_j \wedge \chi_{n-j-1}, j=0,\ldots,n-1$ are all proportional. From \eqref{eq_differential_chi} we obtain that 
\begin{displaymath}
d(\gamma \wedge \chi_0 \wedge \chi_{n-1}) = 0. 
\end{displaymath}
Therefore $\gamma \wedge \chi_0 \wedge \chi_{n-1}$ represents the zero valuation. It follows that each valuation of weight $2$ can be represented by a multiple of the form $\beta \wedge \chi_0 \wedge \chi_{n-1}$, which implies that $\phi_2$ spans $\Val^{SU(n),2}$. 

\subsection{Construction of $\phi_1$}

Suppose that $n=2m$ is even. We claim that 
\begin{displaymath}
\phi_1(K):=\frac{1}{n\pi^m } \nc(K)(\chi_0 \wedge \theta_2^m)
\end{displaymath}
satisfies \eqref{eq_klain_phi1}. 

Since $\nc(P)$ is Legendrian, we can replace $\theta_2$ by $\theta_2+\alpha \wedge \beta=\frac{i}{2} \sum_{j=1}^n dz_j \wedge \bar dz_j$. 

It easily follows that $\phi_1$ is a constant coefficient valuation.

Let $W \in \Gr_n(\mathbb{C}^n)$. If the restriction of the symplectic form to $W$ is non-degenerated, we fix the orientation given by the symplectic form. 

Then 
\begin{displaymath}
[[B_W]]\left((\theta_2+\alpha \wedge \beta)^m\right)=m! \vol(B_W) \prod_{j=1}^m \cos \theta_j(W). 
\end{displaymath}

It thus follows that 
 
\begin{align*}
\phi_1(B_W) & = \frac{1}{n\pi^m} \nc(B_W)(\chi_0 \wedge (\theta_2+\alpha \wedge \beta)^m)\\
& =\frac{1}{n \pi^m} [[B_W]]((\theta_2+\alpha \wedge \beta)^m) [[\partial B_W]]\left(\sum_{j=1}^n (-1)^{j-1} \zeta_j d\zeta_1 \wedge \ldots \wedge \hat{d\zeta_j} \wedge \ldots \wedge d\zeta_n\right) \\
& = \vol(B_W) \prod_{j=1}^m \cos \theta_j(W) \Theta(W). 
\end{align*}

\subsection{Classification of invariant valuations of weight $1$ if $n$ is even} 

Let us next show that $\phi_1$ generates $\Val^{SU(n),1}$. If $\mu \in \Val^{SU(n),1}_k$ with $k \neq n$, then we may apply Klain's injectivity theorem \cite{kl00} (note that $\mu$ is even) and Proposition \ref{prop_orbitspace} (4) to deduce that $\mu$ is $U(n)$-invariant. Therefore $\mu \in \Val^{U(n)} \cap \Val^{SU(n),1}=\{0\}$. 

Hence we may suppose that $\mu \in \Val^{SU(n),1}_n$. Let $\mu$ be represented by a form $\omega$ of bidegree $(n,n-1)$. Proposition \ref{prop_chi_theta} implies that, up to multiples of $\alpha$ and $d\alpha$, $\omega$ is a linear combination of the forms 
\begin{displaymath}
\chi_m \wedge \theta_1^m \quad \text{ and } \quad \beta \wedge \gamma \wedge \chi_{m} \wedge \theta_1^{m-1}.   
\end{displaymath}

We will make frequent use of the relation 
\begin{equation} \label{eq_symplectic_form_relation}
 \theta_s \equiv -\beta \wedge \gamma \mod (\alpha,d\alpha).
\end{equation}

Using \eqref{eq_relation_chi_theta1} and \eqref{eq_symplectic_form_relation}, we obtain that  

\begin{displaymath}
\chi_m \wedge \theta_1^m \equiv m i \beta \wedge \gamma \wedge \chi_m \wedge \theta_1^{m-1} \mod (\alpha,d\alpha).
\end{displaymath}

Hence $\chi_m \wedge \theta_1^m$ and $\beta \wedge \gamma \wedge \chi_m \wedge \theta_1^{m-1}$ induce the same valuation (up to a constant). Since $\phi_1$ is a non-zero valuation of degree $n$ and weight $1$, it follows that $\dim \Val^{SU(n),1}_n=1$. 

\subsection{Classification of invariant valuations of weight $1$ if $n$ is odd} 

Suppose that $n=2m+1$ is odd. Let $\mu$ be a valuation of weight $1$ and of degree $n$. We represent $\mu$ by integration of an invariant differential form $\omega$ of bidegree $(n,n-1)$ and of weight $1$.  

Using Proposition \ref{prop_chi_theta}, we get that $\omega$ is- up to multiples of $\alpha$ and $d\alpha$- a linear combination of the forms  

\begin{displaymath}
\beta \wedge \chi_m \wedge \theta_1^m \quad \text{ and } \quad \gamma \wedge \chi_{m+1} \wedge \theta_1^m.
\end{displaymath}

From \eqref{eq_relation_chi_theta1} we see that $\chi_{m+1} \wedge \theta_1^m$ is divisible by $\theta_s$. Therefore, 
\begin{equation} \label{eq_gamma_term_vanishes}
\gamma \wedge \chi_{m+1} \wedge \theta_1^m \equiv 0 \mod (\alpha,d\alpha).
\end{equation}
 
Using Proposition \ref{prop_park_classification}, \eqref{eq_differential_chi} and \eqref{eq_gamma_term_vanishes}, we compute that 
\begin{displaymath}
d(\chi_{m+1} \wedge \theta_1^m)  \equiv m i \beta \wedge \chi_m \wedge \theta_1^m \mod (\alpha,d\alpha). 
\end{displaymath}

Since multiples of $\alpha$ and $d\alpha$ and exact forms induce the zero valuation, the same holds true for $\beta \wedge \chi_m \wedge \theta_1^m$ and $\gamma \wedge \chi_{m+1} \wedge \theta_1^m$. It follows that $\dim \Val^{SU(n),1}_n=0$. From the Hard Lefschetz Theorem \cite{bebr07} we deduce that $\dim \Val^{SU(n),1}_k=0$ for all $k=0,\ldots,2n$. 

This finishes the proof of Theorem \ref{thm_hadwiger}. 
%----------------------------------------------------------------

\section{Kinematic formulas}
\label{sec_kinematic}

\subsection{Some facts about the product structure and kinematic formulas}

Alesker has shown that there is a dense subspace in the space of all translation invariant valuations on which a natural product structure exists. In fact, this product can even be extended to the much larger space of smooth valuations on a smooth manifold. We refer to \cite{ale04a} and \cite{ale06} for the definition and the properties of this product. 

If $G$ is a compact subgroup of the orthogonal group acting transitively on the unit sphere, then the kinematic formulas for $G$ can be obtained from the product structure of the space $\Val^G$ of $G$-invariant and translation invariant valuations. This important fact is explained in \cite{fu06} and \cite{befu06} and used in a crucial way in the determination of $k_{U(n)}$ in \cite{befu08}. 

Recall that $S\mathbb{C}^n$ is a contact manifold of dimension $2n-1$ with a global contact form $\alpha$. Given an $n-1$-form $\omega$, there exists a unique vertical form $\alpha \wedge \xi$ such that $d(\omega+\alpha \wedge \xi)$ is vertical. The operator $D\omega:=d(\omega+\alpha \wedge \xi)$ is a second-order differential operator which was introduced by Rumin \cite{rum94}. It was first used in integral geometry in \cite{bebr07}. 

In order to compute the product structure on $\Val^{SU(n)}$, we need the following corollary of Theorem 4.1. from \cite{be06}.   

\begin{Proposition} \label{prop_product_formula}
Let $\mu_1,\mu_2 \in \Val^{SU(n)}$ be of degree $n$. Suppose that $\mu_1,\mu_2$ are represented by invariant forms $\omega_1, \omega_2$ respectively. Then 
\begin{displaymath}
\mu_1 \cdot \mu_2 = (-1)^n \frac{2\pi^n}{(n-1)!} c \vol,
\end{displaymath}
where the constant $c$ is determined by 
\begin{displaymath}
\omega_1 \wedge D\omega_2 = c d\vol_{S\mathbb{C}^n}. 
\end{displaymath}
\end{Proposition}

The factor $(-1)^n$ is due to the fact that each $SU(n)$-invariant valuation of degree $n$ is even and thus lies in the $(-1)^n$-eigenspace of the Euler-Verdier involution $\sigma$, compare \cite{be07}. The factor $\frac{2\pi^n}{(n-1)!}$ is the volume of the $2n-1$-dimensional unit sphere. If $\pi:S\mathbb{C}^n \to \mathbb{C}^n$ is the natural projection map, then we have $\pi_* d\vol_{S\mathbb{C}^n}=\frac{2\pi^n}{(n-1)!} d\vol_{\mathbb{C}^n}$.

By Proposition \ref{prop_splitting}, the product of a unitarily invariant valuation of positive degree with $\phi_1$ or $\phi_2$ or their complex conjugates is zero.  
Using Proposition \ref{prop_splitting} and the fact that valuations outside the middle degree are of weight $1$, we obtain that the only non-trivial products are those in weight $0$ as well as $\phi_1 \cdot \bar \phi_1$ and $\phi_2 \cdot \bar \phi_2$. 

\subsection{Computation of $\phi_1 \cdot \bar \phi_1$} 
\label{subs_phi_1} 

Let $n=2m$ be even. We know that $\phi_1$ is represented by $\omega:=\frac{1}{n\pi^m} \chi_0 \wedge \theta_2^m$. 

Let us compute the Rumin differential of $\omega$. Computing modulo $\alpha$ and using \eqref{eq_symplectic_form_relation}, we obtain  
\begin{align*}
d \omega & \equiv \frac{i}{\pi^m} \gamma \wedge \chi_0 \wedge \theta_2^m + \frac{1}{2\pi^m} \beta \wedge \chi_0 \wedge \theta_s \wedge \theta_2^{m-1}\\
& = \frac{i^{m+1}}{(2\pi)^m} \gamma \wedge \chi_m \wedge ((i\theta_1+\theta_s)^m - (i\theta_1-\theta_s)^m) - \frac{1}{2\pi^m} d\alpha \wedge \beta \wedge \chi_0 \wedge \theta_2^{m-1}\\
& = - d\alpha \wedge \left(\frac{2 i^{m+1}}{(2\pi)^m} \gamma \wedge \chi_m \wedge \sum_{j \equiv 1 (2)} i^{m-j}\binom{m}{j}\theta_s^{j-1} \wedge \theta_1^{m-j}+\frac{1}{2\pi^m}\beta \wedge \chi_0 \wedge \theta_2^{m-1}\right).
\end{align*}

Setting
\begin{displaymath}
\xi:=\frac{2 i^{m+1}}{(2\pi)^m} \gamma \wedge \chi_m \wedge \sum_{j \equiv 1 (2)} i^{m-j}\binom{m}{j}\theta_s^{j-1} \wedge \theta_1^{m-j}+\frac{1}{2\pi^m}\beta \wedge \chi_0 \wedge \theta_2^{m-1}
\end{displaymath}
we thus have 
\begin{displaymath}
D \omega=d(\omega+\alpha \wedge \xi)=\alpha \wedge \left( \frac{1}{2\pi^m} \chi_0 \wedge \theta_1 \wedge \theta_2^{m-1}-d\xi\right).
\end{displaymath}

Since we want to compute $\bar \omega \wedge D\omega$, we only need to look at terms in $D\omega$ which are divisible by $\alpha \wedge \beta \wedge \gamma$ and which are not annihilated by $\bar \chi_0$. Since $\bar \chi_0 \wedge \theta_0=0$, there are no terms coming from the differentiation of $\theta_1$ and $\theta_s$ in the first summand of $\xi$. We also note that $\bar \chi_0 \wedge (\theta_s-i\theta_1)=0$. Using the relations from Proposition \ref{prop_chi_theta}, we compute 

\begin{align*}
\bar \omega \wedge D\omega & = \bar \omega \wedge \alpha \wedge \left( \frac{1}{2\pi^m} \chi_0 \wedge \theta_1 \wedge \theta_2^{m-1}-d\xi\right)\\
& = - \bar \omega \wedge \alpha \wedge d\xi\\
& = \frac{2 i^{m+1}}{(2\pi)^m} \bar \omega \wedge \alpha \wedge \gamma \wedge d\chi_m \wedge \sum_{j \equiv 1 (2)} i^{m-j}\binom{m}{j}\theta_s^{j-1} \wedge \theta_1^{m-j}\\
& \quad  +\frac{1}{2\pi^m}  \bar \omega \wedge \alpha \wedge  \beta \wedge d\chi_0 \wedge \theta_2^{m-1}\\
& = \frac{(-1)^m}{(\pi)^m} \bar \omega \wedge \alpha \wedge \gamma \wedge d\chi_m \wedge \theta_1^{m-1}+\frac{mi}{\pi^m}  \bar \omega \wedge \alpha \wedge  \beta \wedge \gamma \wedge \chi_0 \wedge \theta_2^{m-1}\\
& = \frac{(-1)^m mi}{\pi^m} \bar \omega \wedge \alpha \wedge \gamma \wedge \beta \wedge \chi_{m-1} \wedge \theta_1^{m-1}+\frac{mi}{\pi^m}  \bar \omega \wedge \alpha \wedge  \beta \wedge \gamma \wedge \chi_0 \wedge \theta_2^{m-1}\\
& =  \frac{(-1)^m mi}{(2i)^{m-1}\pi^m} \bar \omega \wedge \alpha \wedge \gamma \wedge \beta \wedge \chi_{m-1} \wedge (\theta_s+i\theta_1)^{m-1}\\
& \quad +\frac{mi}{\pi^m}  \bar \omega \wedge \alpha \wedge  \beta \wedge \gamma \wedge \chi_0 \wedge \theta_2^{m-1}\\
& =  \frac{i}{\pi^n} \alpha \wedge \beta \wedge \theta_2^{n-1} \wedge \gamma \wedge \chi_0 \wedge \bar \chi_0.
\end{align*}

One easily sees that $i \gamma \wedge \chi_0 \wedge \bar \chi_0$ equals $2^{n-1}$ times the volume form of the unit sphere in $\mathbb{C}^n$, while $\alpha \wedge \beta \wedge \theta_2^{n-1}$ equals $(n-1)!$ times the volume form of $\mathbb{C}^n$.  

Therefore 
\begin{displaymath}
\bar \omega \wedge D\omega=\frac{2^{n-1}(n-1)!}{\pi^n} d\vol_{S\mathbb{C}^n}. 
\end{displaymath} 
From Proposition \ref{prop_product_formula} we deduce that 
\begin{displaymath}
\bar \phi_1 \cdot \phi_1=2^n \vol.  
\end{displaymath} 

\subsection{Computation of $\phi_2 \cdot \bar \phi_2$}
\label{subs_phi_2}

Let $n \geq 2$ be arbitrary. The valuation $\phi_2$ is represented by the form $\omega:=\frac{(-1)^{n+1}i}{n \omega_n} \beta \wedge \chi_0 \wedge \chi_{n-1}$. We get

\begin{align*}
D(\beta \wedge \chi_0 \wedge \chi_{n-1}) & = d(\beta \wedge \chi_0 \wedge \chi_{n-1} -(n+1)i \alpha \wedge \chi_0 \wedge \chi_{n-1})\\
& =-n(n+2)\alpha \wedge \gamma \wedge \chi_0 \wedge \chi_{n-1},
\end{align*} 

from which we deduce that 
\begin{align*}
\bar \omega \wedge D\omega & = - \frac{n+2}{n \omega_n^2} \beta \wedge \bar \chi_0 \wedge \bar \chi_{n-1} \wedge \alpha \wedge \gamma \wedge \chi_0 \wedge \chi_{n-1}\\
& \quad = \frac{(-1)^{n+1}(n+2)}{n \omega_n^2} \alpha \wedge \beta \wedge \chi_{n-1} \wedge \bar \chi_{n-1} \wedge \gamma \wedge  \chi_0 \wedge \bar \chi_{0}.
\end{align*}

The form $\alpha \wedge \beta \wedge \chi_{n-1} \wedge \bar \chi_{n-1}$ is $2^{n-1}i^{n^2-1}$ times the volume form of $\mathbb{C}^n$; while $\gamma \wedge  \chi_0 \wedge \bar \chi_{0}$ equals $2^{n-1}i^{n^2-1}$ times the volume form of the unit sphere $S^{2n-1}$. Hence 
\begin{displaymath}
\bar \omega \wedge D\omega = \frac{n+2}{n \omega_n^2} 4^{n-1} d\vol_{S\mathbb{C}^n} 
\end{displaymath}
and thus 
\begin{displaymath}
\bar \phi_2 \cdot \phi_2 = \frac{(-1)^{n}(n+2) 2^{2n-1} \omega_{2n}}{\omega_n^2}  \vol.  
\end{displaymath}

We remark that the pairing $(\mu_1,\mu_2) \mapsto \bar \mu_1 \cdot \mu_2$ on $\Val_n^{SU(n)}$ is not positive definite if $n$ is odd. Equivalently, the pairing $(\mu_1,\mu_2) \mapsto \mu_1 \cdot \mu_2$ on {\it real-valued} valuations in $\Val_n^{SU(n)}$ is not positive definite. In contrast to this, it was shown in \cite{befu08} that the restriction of this pairing to $\Val^{U(n)}_n$ is positive definite for all $n$.  

\subsection{Kinematic formulas}
The principal kinematic formula for $SU(n)$ follows from \ref{subs_phi_1} and \ref{subs_phi_2} and Theorem 2.6. of \cite{fu06}:

\begin{align*}
 k_{SU(n)}(\chi)& =k_{U(n)}(\chi) +\frac{(-1)^{n}\omega_n^2}{(n+2) 2^{2n-1} \omega_{2n}}  (\phi_2 \otimes \bar \phi_2+\bar \phi_2 \otimes \phi_2) \\
& \quad \left(+\frac{1}{2^n}(\phi_1 \otimes \bar \phi_1+\bar \phi_1 \otimes \phi_1) \text{ if } n \text{ is even}\right).
\end{align*}

If $\mu$ is a $U(n)$-invariant valuation of degree $k>0$, then $\mu \cdot \phi_1=0$ and $\mu \cdot \phi_2=0$ by Proposition \ref{prop_splitting} and the fact that valuations of non-zero weight appear only in degree $n$. Using Lemma 2.4. of \cite{fu06}, we obtain 
\begin{displaymath}
 k_{SU(n)}(\mu)=k_{U(n)}(\mu). 
\end{displaymath}

By the same lemma, 
\begin{displaymath}
 k_{SU(n)}(\phi_2)=\phi_2 \otimes \vol + \vol \otimes \phi_2; 
\end{displaymath}
for all $n$, and 
\begin{displaymath}
 k_{SU(n)}(\phi_1)=\phi_1 \otimes \vol + \vol \otimes \phi_1 
\end{displaymath}
for even $n$. 

This establishes the whole array of kinematic formulas for $SU(n)$ and finishes the proof of Theorem \ref{thm_kinematic}. 

\subsection{Additive kinematic formula}

\begin{Proposition}
For all compact convex sets $K,L \subset \mathbb{C}^n$, the following additive kinematic formula holds. 
\begin{align*}
\int_{SU(n)} \vol(K + g L) dg & =\int_{U(n)} \vol(K + gL)dg \\
& \quad + \frac{(-1)^{n}\omega_n^2}{(n+2) 2^{2n-1} \omega_{2n}}  (\phi_2(K)\bar \phi_2(L)+\bar \phi_2(K) \phi_2(L)) \\
& \quad \left(+\frac{1}{2^n}(\phi_1(K) \bar \phi_1(L)+\bar \phi_1(K) \phi_1(L)) \text{ if } n \text{ is even}\right).
\end{align*}
\end{Proposition}

\proof
Since the Fourier transform acts trivially on $\Val_n^{SU(n)}$, the proposition follows immediately from Theorem 1.7. of \cite{befu06}.
\endproof

Higher additive kinematic formulas (where $\vol$ under the integral is replaced by another $SU(n)$-invariant valuations) can be obtained from the proposition by applying the results of \cite{befu06}. We leave the details to the reader. 
%------------------------------------------

\end{document}